\documentclass[11 pt]{article}
\usepackage{amsmath}
\usepackage{amsfonts}
\usepackage{color}

\begin{document}

\vspace*{.5cm}

\begin{center}
{\large \textbf{Lightlike Submanifolds of Metallic Semi-Riemannian Manifolds}%
}

\bigskip

{\small Bilal Eftal ACET, Feyza Esra ERDO\u{G}AN and Selcen Y\"{U}KSEL PERKTA\c{S}}
\end{center}

\noindent \textbf{Abstract}. {\small Our aim in this paper is to investigate
some types of lightlike submanifolds in metallic semi-Riemannian manifolds.
We study invariant and screen semi-invariant lightlike submanifolds of
metallic semi-Riemannian manifolds and give examples. We obtain the
conditions for the induced connection to be a metric connection. Also, we
find necessary and sufficient conditions for the distributions involved in
the definitions of such submanifolds to be integrable.}

\medskip


\section*{\textbf{1.~Introduction}}

\setcounter{equation}{0} \renewcommand{\theequation}{1.\arabic{equation}}

The existence of a wide range of applications in mathematics and physics has
made the Riemannian and semi-Riemannian geometries an important research
area for differential geometry. But by the 1970s, Einstein's general
relativity theory shifted the interest on Riemannian and semi-Riemannian
geometries to the Lorentzian geometry. Later, Duggal and Bejancu published a
book on lightlike geometry in 1996 \cite{Duggal-Bejancu} and filled an
important gap in the theory of submanifolds. In this book, the geometric
objects for obtaining the Gauss-Codazzi equations of a lightlike submanifold
are defined by means of a non-degenerate screen distribution.

The main difference between the theory of lightlike submanifolds and
semi-Riemannian submanifolds arises due to the fact that in the first case,
a part of the normal vector bundle $TN^{\perp }$ lies in the tangent bundle $%
TN$ of the submanifold $N$ of a semi-Riemannian manifold $\breve{N},$
whereas in the second case $TN\cap TN^{\perp }=\{0\}.$ Thus, the basic
problem of lightlike submanifolds is to replace the intersecting part by a
vector subbundle whose sections are nowhere tangent to $N.$ To construct a
nonintersecting lightlike transversal vector bundle of the tangent bundle,
Duggal and Bejancu used an extrinsic approach while K\"{u}peli used an
intrinsic approach \cite{Kupeli}. Since then, many authors have studied the
geometry of lightlike hypersurfaces and lightlike submanifolds. Recent
studies have been updated in \cite{Duggal-Sahin-2}. Many studies on \ lightlike submanifolds
have been reported by many geometers (see \cite{APK3, APK1, FBR-1, FC, Burcin-1} and the references therein). In this paper, we follow
the approach given by Duggal and Bejancu in \cite{Duggal-Bejancu}. We note
that lightlike hypersurfaces are examples of physical models of Killing
horizons in general relativity \cite{Galloway} and the relationship between
Killing horizons and black holes is based on Hawking's area theorem \cite%
{Hawking}.

It is known that the number $\phi =(1+\sqrt{5})/2=1,618033...$ is a solution
of the equation $x^{2}-x-1=0$ and it is called as golden ratio. The golden
ratio is very interesting because of its use in art works and frequent
occurrence in the nature. Thus, Crasmareanu and Hretcanu defined a Golden
manifold $N$ by a tensor field $\Phi $ on $N$ satisfies $\Phi ^{2}=\Phi +I$
in \cite{Crasmareanu-Hretcanu}. In the same paper, the authors showed that $%
\phi $ and $1-\phi $ are eigenvalues of $\Phi .$ Then, in \cite{Sahin-Akyol} 
\c{S}ahin and Akyol introduced Golden maps between Golden Riemannian
manifolds and showed that such maps are harmonic maps. Finally, lightlike
hypersurfaces of a Golden semi-Riemannian manifolds was studied by Poyraz
and Ya\c{s}ar in \cite{Poyraz-Yasar}.

In 1997 Spinadel introduced metallic means family or metallic proportions as
a generalization of the golden mean in \cite{Spinadel}-\cite{Spinadel-4}.
Let $p$ and $q$ be positive integers. Then, member of the metallic means
family is the positive solution of the equation $x^{2}-px-q=0$ and these
numbers, which are called $(p,q)$ metallic numbers, denoted by 
\begin{equation}
\sigma _{p,q}=\frac{p+\sqrt{p^{2}+4q}}{2}.  \label{eq:1.1}
\end{equation}%
The members of the metallic means family take the name of a metal, like the
golden mean, the silver mean, the bronze mean, the copper mean and many
others. A metallic manifold $\breve{N}$ has a tensor field $\breve{J}$ such
that the equality $\breve{J}^{2}=p\breve{J}+qI$ is satisfied, where $p$ and $%
q$ positive integers and the eigenvalues of automorphism $\breve{J}$ of the
tangent bundle $TN$ are $\sigma _{p,q}$ and $p-\sigma _{p,q}$ \cite%
{Crasmareanu-Hretcanu-2}. Moreover, if $(\breve{N},\breve{g}%
)$ is a Riemannian manifold endowed with a metallic structure $\breve{J}$
such that the Riemannian metric $\breve{g}$ is $\breve{J}$-compatible, i.e., 
$\breve{g}(\breve{J}X,Y)=\breve{g}(X,\breve{J}Y),$ for any $X,Y\in \chi (%
\breve{N})$, then $(\breve{g},\breve{J})$ is called metallic Riemannian
structure and $(\breve{N},\breve{g},\breve{J})$ is a metallic Riemannian
manifold. Metallic structure on the ambient Riemannian manifold provides
important geometrical results on the submanifolds, since it is an important
tool while investigating the geometry of submanifolds. Invariant,
anti-invariant, semi-invariant, slant and semi-slant submanifolds of a
metallic Riemannian manifold are studied in \cite{Adara-1, Adara-2, Adara-3}%
. Some types of lightlike submanifolds of a metallic semi-Riemannian
manifold are introduced in \cite{BEA, ER} and the authors obtained important
characterizations on such submanifolds of metallic semi-Riemannian manifolds
with examples.

In this paper, we study some special types of lightlike submanifolds in
metallic semi-Riemannian manifolds. In section~2, we give basic informations
needed for the rest of the paper. In section~3 and section~4, we research
invariant lightlike submanifolds and screen semi-invariant lightlike
submanifolds of metallic semi-Riemannian manifolds, respectively. In these
sections, we give some characterizations and investigate the geometry of
leaves of distributions which arise from definitions. In general, the
induced connection of a lightlike submanifold is not a metric connection.
Therefore it is an important problem to find conditions for the induced
connection to be a metric connection. So, we find necessary and sufficient
conditions for the induced connection to be a metric connection. Finally,
note that the paper contains examples.

\section*{\textbf{2.\, Preliminaries}}

\setcounter{equation}{0} \renewcommand{\theequation}{2.\arabic{equation}}

Let $\breve{N}$ be a differentiable manifold and $\breve{J}$ be a $(1,1)$
type tensor field on $\breve{N}.$ If the following equation is satisfied,
then $\breve{J}$ is called a metallic structure on $\breve{N}:$ 
\begin{equation}
\breve{J}^{2}=p\breve{J}+qI,  \label{eq:2.1a}
\end{equation}%
where $p,q$ are positive integers and $I$ is the identity operator on the
Lie algebra $\chi (\breve{N})$ of the vector fields on $\breve{N}.$ If $%
\breve{J}$ is a self-adjoint operator with respect to semi-Riemann metric $%
\breve{g}$ of a semi-Riemann manifold $\breve{N},$ that is, 
\begin{equation}
\breve{g}(\breve{J}U,V)=\breve{g}(U,\breve{J}V),  \label{eq:2.1b}
\end{equation}%
is satisfied, then $\breve{g}$ is said to be $\breve{J}$-compatible and $(%
\breve{N},\breve{J},\breve{g})$ is called a metallic semi-Riemannian
manifold. Using (\ref{eq:2.1b}), we can write 
\begin{equation}
\breve{g}(\breve{J}U,\breve{J}V)=p\breve{g}(U,\breve{J}V)+q\breve{g}(U,V),
\label{eq:2.1c}
\end{equation}%
for any $U,V\in \Gamma (T\breve{N}).$

It is well known that the non-degenerate metric $\breve{g}$ of a $(m+n)$%
-dimensional semi-Riemann manifold $\breve{N}$ is not always induced as a
non-degenerate metric on an $m$-dimensional submanifold $N$ of $\breve{N}.$
If the induced metric $g$ is degenerate on $N$ and $rank(Rad(TN))=r,$ $1\leq
r\leq m,$ then $(N,g)$ is called a lightlike submanifold of $(\breve{N},%
\breve{g}),$ where the radical distribution $Rad(TN)$ and the normal bundle $%
TN^{\perp }$ of the tangent bundle $TN$ are defined by 
\begin{equation*}
Rad(TN)=TN\cap TN^{\perp }
\end{equation*}%
and 
\begin{equation*}
TN^{\perp }=\cup _{x\in N}\{u\in T_{x}\,\breve{N}\mid \breve{g}(u,v)=0,\,\,\
\forall v\in T_{x}\,N\}.
\end{equation*}%
Since $TN$ and $TN^{\perp }$ are degenerate vector subbundles, there exist
complementary non-degenerate distributions $S(TN)$ and $S(TN^{\bot })$ of $%
Rad(TN)$ in $TN$ and $TN^{\perp },$ respectively, which are called the
screen distribution and screen transversal bundle (or co-screen
distribution) of $N$ such that 
\begin{equation*}
TN=S(TN)\perp Rad(TN)\,,\,\,\ TN^{\bot }=S(TN^{\bot })\perp Rad(TN).
\end{equation*}%
On the other hand, consider an orthogonal complementary bundle $S(TN)^{\bot
} $ to $S(TN)$ in $T\breve{N}$ such that 
\begin{equation}
S(TN)^{\bot }=S(TN^{\bot })\bot S(TN^{\bot })^{\bot },  \label{eq:2.a}
\end{equation}%
where $S(TN^{\bot })^{\bot }$ is the orthogonal complementary to $S(TN^{\bot
})$ in $S(TN)^{\bot }.$

We now recall the following important result. \newline

\noindent \textbf{Theorem~2.1.~}\textsl{Let }$(N,g,S(TN),S(TN^{\bot }))$%
\textsl{\ be a $r$-lightlike submanifold of a semi-Riemannian manifold $(%
\breve{N},\breve{g}).$ Then, there exists a complementary vector bundle $%
ltr(TN)$ called a lightlike transversal bundle of $Rad(TN)$ in $S(TN^{\bot
})^{\bot }$ and a basis of $\Gamma (ltr(TN)\mid _{U})$ consists of smooth
sections $\{N_{1},...,N_{r}\}$ of $S(TN^{\bot })^{\bot }\mid _{U}$ such that 
\begin{equation*}
\breve{g}(\xi _{i},N_{j})=\delta _{ij}\,,\,\,\ \breve{g}(N_{i},N_{j})=0\,,\,%
\,\ i,j=1,..,r,
\end{equation*}%
where $\{\xi _{1},...,\xi _{r}\}$ is a basis of $\Gamma (Rad(TN))$} \cite[%
page 144]{Duggal-Bejancu}. \newline

This result implies that there exists a complementary (but not orthogonal)
vector bundle $tr(TN)$ to $TN$ in $T\breve{N}|_{N},$ which is called
transversal vector bundle, such that the following decompositions are hold: 
\begin{equation}
tr(TN)=ltr(TN)\perp S(TN^{\perp })  \label{eq:2.b}
\end{equation}%
and 
\begin{equation}
S(TN^{\bot })^{\bot }=Rad(TN)\oplus ltr(TN).  \label{eq:2.c}
\end{equation}%
Thus, using (\ref{eq:2.a}), (\ref{eq:2.b}) and (\ref{eq:2.c}), we get 
\begin{align}
T\breve{N}|_{N}& =S(TN)\perp S(TN)^{\bot }  \notag \\
& =S(TN)\perp \left\{ Rad(TN)\oplus ltr(TN)\right\} \perp S(TN^{\perp }) 
\notag \\
& =TN\oplus tr(TN).  \label{eq:2.d}
\end{align}

A submanifold $(N,g,S(TN),S(TN^{\perp })$ is called \newline
(1):\thinspace\ $r$\thinspace -\thinspace lightlike if $r\,<\min \{m,\,n\},$ 
\newline
(2):\thinspace\ Co\thinspace -\thinspace isotropic if $r\,=\,n\,<\,m,\,\,\ $%
i.e.,$\,\,\ S(TN^{\perp })=\{0\},$\newline
(3):\thinspace\ Isotropic if $r\,=\,m\,<\,n,\,\,\ $i.e.,$\,\,\ S(TN)=\{0\}$
and \newline
(4):\thinspace\ Totally lightlike if $r\,=\,m\,=\,n,\,\,\ $i.e.,$\,\,\
S(TN)=\{0\}=S(TN^{\perp }).$ \newline

The Gauss and Weingarten equations of $N$ are given by 
\begin{equation}
\bar{\nabla}_{U}V=\nabla _{U}V+h(U,V)\,\,,\,\,\ \forall U,V\in \Gamma (TN),
\label{eq:2.4}
\end{equation}%
and 
\begin{equation}
\bar{\nabla}_{U}N=-A_{N}U+\nabla _{U}^{t}N\,\,\,,\,\,\ \forall U\in \Gamma
(TN)\,\,,\,\ N\in \Gamma (tr(TN)),  \label{eq:2.5}
\end{equation}%
where $\{\nabla _{U}V,A_{N}U\}$ and $\{h(U,V),\nabla _{U}^{t}N\}$ are belong
to $\Gamma (TN)$ and $\Gamma (tr(TN)),$ respectively. $\nabla $ and $\nabla
^{t}$ are linear connections on $N$ and on the vector bundle $tr(TN),$
respectively. The second fundamental form $h$ is a symmetric $\mathcal{F}(N)$%
-bilinear form on $\Gamma (TN)$ with values in $\Gamma (tr(TN))$ and the
shape operator $A_{V}$ is a linear endomorphism of $\Gamma (TN).$

If we consider (\ref{eq:2.d}) and using the projection morphisms denoted by 
\begin{equation*}
L:tr(TN)\rightarrow ltr(TN)\,,\,\,\,\,\ S:tr(TN)\rightarrow S(TN^{\perp }),
\end{equation*}%
then, for any $U,V\in \Gamma (TN),$ $N\in \Gamma (ltr(TN))$ and $W\in \Gamma
(S(TN^{\perp })),$ we can write 
\begin{eqnarray}
\bar{\nabla}_{U}V &=&\nabla _{U}V+h^{\l }(U,V)+h^{s}(U,V),  \label{eq:2.6} \\
\bar{\nabla}_{U}N &=&-A_{N}U+\nabla _{U}^{\l }(N)+D^{s}(U,N),  \label{eq:2.7}
\\
\bar{\nabla}_{U}W &=&-A_{W}U+D^{l}(U,W)+\nabla _{U}^{s}(W),  \label{eq:2.8}
\end{eqnarray}%
where $\{\nabla _{U}^{\l }N,D^{l}(U,W)\}$ and $\{D^{s}(U,N),\nabla
_{U}^{s}W\}$ are parts of $ltr(TN)$ and $S(TN^{\bot }),$ respectively and $%
h^{\l }(U,V)=Lh(U,V)\in \Gamma (ltr(TN)),$ $h^{s}(U,V)=Sh(U,V)\in \Gamma
(S(TN^{\bot }))$. Denote the projection of $TN$ on $S(TN)$ by $\bar{P}.$
Then, by using (\ref{eq:2.4}), (\ref{eq:2.6})-(\ref{eq:2.8}) and taking
account that $\bar{\nabla}$ is a metric connection we obtain 
\begin{equation}
\breve{g}(h^{s}(U,V),W)+\breve{g}(Y,D^{l}(U,W))=g(A_{W}U,Y),  \label{eq:2.9}
\end{equation}%
\begin{equation}
\breve{g}(D^{s}(U,N),W)=\breve{g}(N,A_{W}U),  \label{eq:2.10}
\end{equation}%
\begin{equation}
{\nabla }_{U}{\bar{P}}V={\nabla }_{U}^{\ast }{\bar{P}}V+h^{\ast }(U,{\bar{P}}%
V),  \label{eq:2.11}
\end{equation}%
and 
\begin{equation}
{\nabla }_{U}{\xi }=-A_{\xi }^{\ast }U+{{\nabla }^{\ast }}_{U}^{t}{\xi ,}
\label{eq:2.12}
\end{equation}%
for $U,V\in \Gamma (TN)$ and $\xi \in \Gamma (RadTN),$ where ${\nabla }%
^{\ast }$ and ${{\nabla }^{\ast }}^{t}$ are induced connections on $S(TN)$
and $Rad(TN)$, respectively. On the other hand, $h^{\ast }$ and $A^{\ast }$
are $\Gamma (Rad(TN))$-valued and $\Gamma (S(TN))$-valued $\mathcal{F}(N)$%
-bilinear forms on $\Gamma (TN)\times \Gamma (S(TN))$ and $\Gamma
(Rad(TN))\times \Gamma (TN),$ respectively. $h^{\ast }$ is called local
second fundamental form on $S(TN)$ and $A^{\ast }$ is second fundamental
form of $Rad(TN).$

By using above equations we obtain 
\begin{eqnarray}
\breve{g}(h^{l}(U,{\bar{P}}V),\xi ) &=&g(A_{\xi }^{\ast }U,{\bar{P}}V),
\label{eq:2.13} \\
\breve{g}(h^{\ast }(U,{\bar{P}}V),N) &=&g(A_{N}U,{\bar{P}}V),
\label{eq:2.14} \\
\breve{g}(h^{l}(U,\xi ),\xi )=0, &&A_{\xi }^{\ast }{\xi }=0.  \label{eq:2.15}
\end{eqnarray}

In general, the induced connection $\nabla $ on $N$ is not a metric
connection. Since $\bar{\nabla}$ is a metric connection, by using (\ref%
{eq:2.6}), we get 
\begin{equation}
(\nabla _{U}g)(V,Z)=\breve{g}(h^{l}(U,V),Z)+\breve{g}(h^{l}(U,Z),V).
\label{eq:2.16}
\end{equation}%
However, it is important to note that $\nabla ^{\star }$ is a metric
connection on $S(TN).$

\section*{\textbf{3.\, Invariant Lightlike Submanifolds of Metallic
Semi-Riemannian Manifolds}}

\setcounter{equation}{0} \renewcommand{\theequation}{3.\arabic{equation}}

\noindent \textbf{Definition~1.~}\textsl{Let }$(\breve{N},\breve{J},\breve{g})$
\textsl{\ be a metallic semi-Riemannian manifold and $(N,g)$ be
a lightlike submanifold of $\breve{N}.$ Then, we say that $N$ is an
invariant lightlike submanifold of $\breve{N},$ if the following conditions
are satisfied:} 
\begin{equation}
\breve{J}(S(TN))=S(TN)  \label{eq:3.1}
\end{equation}%
and 
\begin{equation}
\breve{J}(Rad(TN))=Rad(TN).  \label{eq:3.2}
\end{equation}

\medskip

\noindent \textbf{Corollary~3.1~} \textsl{Let $(\breve{N},\breve{J},\breve{g}%
)$ be a metallic semi-Riemannian manifold and $(N,g)$ be an invariant
lightlike submanifold of $\breve{N}.$ Then, the lightlike transversal distribution $ltr(TN)$ is
invariant with respect to $\breve{J}.$ } \newline
\noindent \textbf{Proof.~} We assume that $N$ is an invariant lightlike
submanifold of $\breve{N}.$ Then, for any $U\in \Gamma (S(TN)),$ $\xi \in
\Gamma (Rad(TN))$ and $N\in \Gamma (ltr(TN)),$ we have 
\begin{equation}
\breve{g}(\breve{J}N,\xi )=\breve{g}(N,\breve{J}\xi )\neq 0  \label{eq:3.3}
\end{equation}%
and 
\begin{equation}
\breve{g}(\breve{J}N,U)=\breve{g}(N,\breve{J}U)=0.  \label{eq:3.4}
\end{equation}%
From (\ref{eq:3.3}) and (\ref{eq:3.4}), it is clear that $\breve{J}N$ have
no components in $\Gamma (Rad(TN))$ and $\Gamma (S(TN)),$ respectively. On
the other hand, since $Rad(TN)\bot S(TN^{\bot }),$ then $S(TN^{\bot })$ does
not contain $\breve{J}N.$ Thus, the proof is completed. \newline

\medskip

\noindent \textbf{Example~1.~} Let $\breve{N}=R_{1}^{5}$ be a
metallic semi-Riemannian manifold of signature $(-,+,+,+,+)$ with respect to
the canonical basis $\{\frac{\partial }{\partial x_{1}},\frac{\partial }{%
\partial x_{2}},\frac{\partial }{\partial x_{3}},\frac{\partial }{\partial
x_{4}},\frac{\partial }{\partial x_{5}}\}.$ \newline
Consider a metallic structure $\breve{J}$ defined by 
\begin{equation*}
\breve{J}(x_{1},x_{2},x_{3},x_{4},x_{5})=(\sigma x_{1},\sigma x_{2},\sigma
x_{3},(p-\sigma )x_{4},(p-\sigma )x_{5}).
\end{equation*}%
Let $N$ be a submanifold of $(R_{1}^{5},\breve{J},\breve{g})$ given by 
\begin{equation*}
x_{1}=u_{3}\,,\,\,\ x_{2}=-\sin \alpha u_{1}+\cos \alpha u_{3},
\end{equation*}
\begin{equation*}
x_{3}=\cos \alpha u_{1}+\sin \alpha u_{3}\,,\,\,\ x_{4}=u_{2}\,,\,\,\
x_{5}=0.
\end{equation*}%
Then $TN$ is spanned by $\{Z_{1},Z_{2},Z_{3}\}$, where 
\begin{equation*}
Z_{1}=-\sin \alpha {\partial \,x_{2}}+\cos \alpha {\partial \,x_{3}},
\end{equation*}
\begin{equation*}
Z_{2}={\partial \,x_{4}},
\end{equation*}
\begin{equation*}
Z_{3}={\partial \,x_{1}}+\cos\alpha {\partial \,x_{2}}+\sin \alpha {\partial \,x_{3}}.
\end{equation*}%
Hence $N$ is a $1-$lightlike submanifold of $R_{1}^{5}$ with 
$$
Rad(TN)=Span\{Z_{3}\}
$$
and 
$$
S(TN)=Span\{Z_{1},Z_{2}\}.
$$ 
It is easy to see
that 
$$
\breve{J}Z_{3}=\sigma Z_{3}\in \Gamma (Rad(TN)),\,\,\, \breve{J}%
Z_{1}=\sigma Z_{1}\in \Gamma (S(TN)),
$$
$$
\breve{J}Z_{2}=\sigma Z_{2}\in\Gamma (S(TN)),
$$
which mean that $S(TN)$ and $Rad(TN)$ is invariant with
respect to $\breve{J}.$ On the other hand, by direct calculations, we get
the lightlike transversal bundle and screen-transversal distribution are
spanned by 
\begin{equation*}
N=\frac{1}{2}\{-\partial \,x_{1}+\cos \alpha \partial \,x_{2}+\sin \alpha
\partial \,x_{3}\}\,,\,\,\,\ W=\partial \,x_{5},
\end{equation*}%
respectively. It is clear that $ltr(TN)$ and $S(TN^{\bot })$ are invariant
distributions. Thus, $N$ is an invariant lightlike submanifold of $\breve{N}$.\newline

\medskip

Let $(\breve{N},\breve{J},\breve{g})$ be a metallic semi-Riemannian manifold
and $N$ be an invariant lightlike submanifold of $\breve{N}.$ In this paper,
we assume that 
\begin{equation}
\bar{\nabla}\breve{J}=0,  \label{eq:3.5}
\end{equation}%
which implies that $\bar{\nabla}_{U}${$\breve{J}$}${V}=\breve{J}\bar{\nabla}%
_{U}V$ and similarly $\bar{\nabla}_{U}${$\breve{J}$}${N}=\breve{J}\bar{\nabla%
}_{U}N,$ for all $U,V\in \Gamma (TN),$ $N\in \Gamma (ltr(TN)).$ \newline

Now we denote the projection morphisms on $S(TN)$ and $Rad(TN)$ by $T$ and $%
Q,$ respectively. Then, for any $U\in \Gamma (TN)$ we write 
\begin{equation}
U=TU+QU,  \label{eq:3.6}
\end{equation}%
where $TU\in \Gamma (S(TN))$ and $QU\in \Gamma (Rad(TN)).$

Appliying $\breve{J}$ to (\ref{eq:3.6}) we get 
\begin{equation}
\breve{J}U=\breve{J}TU+\breve{J}QU.  \label{eq:3.7}
\end{equation}%
If we denote $\breve{J}TU$ and $\breve{J}QU$ by $SU$ and $LU,$ respectively,
then we can rewrite (\ref{eq:3.7}) as 
\begin{equation}
\breve{J}U=SU+LU,  \label{eq:3.8}
\end{equation}%
where $SU\in \Gamma (S(TN))$ and $LU\in \Gamma (Rad(TN)).$ \newline

Let $N$ be an invariant lightlike submanifold of a metallic semi-Riemannian
manifold $\breve{N}.$ If we differentiate (\ref{eq:3.8}) and using (\ref%
{eq:3.5}), (\ref{eq:2.4}), (\ref{eq:2.6})-(\ref{eq:2.8}), for any $U,V\in
\Gamma (TN),$ we have 
\begin{eqnarray}
S\nabla _{U}V+L\nabla _{U}V+\breve{J}h^{l}(U,V)+\breve{J}h^{s}(U,V)
&=&\nabla _{U}^{\ast }SV+h^{\ast }(U,SV)  \notag \\
&+&h^{l}(U,SV)+h^{s}(U,SV)  \notag \\
&-&A_{LV}^{\ast }U+\nabla _{U}^{\ast t}LV  \label{eq:3.9} \\
&+&h^{l}(U,LV)+h^{s}(U,LV).  \notag
\end{eqnarray}%
Considering the tangential, lightlike transversal and screen transversal
parts of this equation we obtain the following

\medskip

\noindent \textbf{Lemma~3.1} \textsl{Let }$(\breve{N},\breve{J},\breve{g})$%
{\textsl{\ be a metallic semi-Riemannian manifold and }}$(N,g)$%
{\textsl{\ be an invariant lightlike submanifold of }}$\breve{N}$%
{\textsl{$.$ Then, we have} } 
\begin{equation}
S\nabla _{U}V=\nabla _{U}^{\ast }SV-A_{LV}^{\ast }U,  \label{eq:3.10}
\end{equation}%
\begin{equation}
L\nabla _{U}V=h^{\ast }(U,SV)+\nabla _{U}^{\ast t}LV,  \label{eq:3.11}
\end{equation}%
\begin{equation}
\breve{J}h^{l}(U,V)=h^{l}(U,\breve{J}V),  \label{eq:3.12}
\end{equation}%
\begin{equation}
\breve{J}h^{s}(U,V)=h^{s}(U,\breve{J}V),  \label{eq:3.13}
\end{equation}%
\begin{equation}
\breve{J}\nabla _{U}V=\nabla _{U}^{\ast }SV-A_{LV}^{\ast }U+h^{\ast
}(U,SV)+\nabla _{V}^{\ast t}LV,  \label{eq:3.14}
\end{equation}%
\textsl{where} $U,V\in \Gamma (TN)$.

\medskip

\noindent \textbf{Theorem~3.1.~}\textsl{Let }$N$\textsl{\ be an invariant
lightlike submanifold of a metallic semi-Riemannian manifold }$\breve{N}.$%
\textsl{\ Then, the radical distribution }$Rad(TN)$\textsl{\ is integrable
if and only if either }%
\begin{equation*}
A_{\breve{J}X}^{\ast }Y=A_{\breve{J}Y}^{\ast }X\,\,\ and\,\,\ A_{X}^{\ast
}Y=A_{Y}^{\ast }X
\end{equation*}%
\textsl{or }%
\begin{equation*}
A_{\breve{J}X}^{\ast }Y-A_{\breve{J}Y}^{\ast }X=p(A_{X}^{\ast }Y-A_{Y}^{\ast
}X),
\end{equation*}%
\textsl{for any }$X,Y\in \Gamma (Rad(TN))$\textsl{\ and }$Z\in \Gamma
(S(TN)).$

\medskip

\noindent \textbf{Proof.~} We know that the distribution $Rad(TN)$ is
integrable if and only if $[X,Y]\in \Gamma (Rad(TN)),$ for all $X,Y\in
\Gamma (Rad(TN)),$ that is, 
\begin{equation*}
\breve{g}([X,Y],Z)=0.
\end{equation*}%
Thus, for any $Z\in \Gamma (S(TN)),$ using (\ref{eq:2.1c}) and (\ref{eq:3.5}%
), we have 
\begin{eqnarray*}
\breve{g}([X,Y],Z) &=&\breve{g}(\bar{\nabla}_{X}Y-\bar{\nabla}_{Y}X,Z) \\
&=&\frac{1}{q}\breve{g}(\bar{\nabla}_{X}{\breve{J}Y}-\bar{\nabla}_{Y}{\breve{%
J}}X,{\breve{J}Z}) \\
&&+\frac{p}{q}\breve{g}(\bar{\nabla}_{Y}X-\bar{\nabla}_{X}Y,{\breve{J}Z}) \\
&=&0.
\end{eqnarray*}%
If we use (\ref{eq:2.6}), we get 
\begin{eqnarray*}
0=\breve{g}(\nabla _{X}{\breve{J}Y},\breve{J}Z)-p\breve{g}(\nabla _{X}Y,%
\breve{J}Z)-\breve{g}(\nabla _{Y}{\breve{J}}X,\breve{J}Z) 
+p\breve{g}(\nabla _{Y}X,\breve{J}Z).
\end{eqnarray*}%
Finally, using (\ref{eq:2.12}) in the last equation, we obtain 
\begin{equation*}
\breve{g}(-A_{\breve{J}Y}^{\ast }X+A_{\breve{J}X}^{\ast }Y,\breve{J}Z)+p%
\breve{g}(A_{Y}^{\ast }X-A_{X}^{\ast }Y,\breve{J}Z)=0,
\end{equation*}%
which completes the proof. \newline

\medskip

\noindent \textbf{Theorem~3.2.~}\textsl{Let }$N$\textsl{\ be an invariant
lightlike submanifold of a metallic semi-Riemannian manifold }$\breve{N}.$%
\textsl{\ Then, the screen distribution }$S(TN)$\textsl{\ is integrable if
and only if either }$h^{\ast }$\textsl{\ is symmetric and self-adjoint or }%
\begin{equation*}
h^{\ast }(U,\breve{J}V)-h^{\ast }(V,\breve{J}U)=p(h^{\ast }(U,V)-h^{\ast
}(V,U)),
\end{equation*}%
\textsl{for all }$U,V\in \Gamma (S(TN))$\textsl{\ and }$N\in \Gamma
(ltr(TN)).$\textsl{\ }\newline

\medskip

\noindent \textbf{Proof.~} $S(TN)$ is integrable if and only if $\breve{g}%
([U,V],N)=0,$ for all $U,V\in \Gamma (S(TN))$ and $N\in \Gamma (ltr(TN)).$
Using (\ref{eq:2.1c}) and (\ref{eq:3.5}), we get 
\begin{eqnarray*}
\breve{g}([U,V],N) &=&\breve{g}(\bar{\nabla}_{U}V-\bar{\nabla}_{V}U,N) \\
&=&\frac{1}{q}\breve{g}(\bar{\nabla}_{U}{\breve{J}V}-\bar{\nabla}_{V}{\breve{%
J}U},{\breve{J}N}) \\
&&+\frac{p}{q}\breve{g}(\bar{\nabla}_{V}U-\bar{\nabla}_{U}V,{\breve{J}N}) \\
&=&0.
\end{eqnarray*}%
If we use (\ref{eq:2.11}), we have 
\begin{eqnarray*}
&&\breve{g}(h^{\ast }(U,\breve{J}V)-h^{\ast }(V,\breve{J}U),\breve{J}N) \\
&&-p\breve{g}(h^{\ast }(U,V)-h^{\ast }(V,U),\breve{J}N) \\
&=&0,
\end{eqnarray*}%
which completes the proof.

\medskip

\noindent \textbf{Theorem~3.3.~}\textsl{Let }$N$\textsl{\ be an invariant
lightlike submanifold of a metallic semi-Riemannian manifold }$\breve{N}.$%
\textsl{\ Then, the induced connection }$\nabla $\textsl{\ on }$N$\textsl{\
is a metric connection if and only if }%
\begin{equation*}
A_{\breve{J}\xi }^{\ast }U=pA_{\xi }^{\ast }U,
\end{equation*}%
\textsl{for \ all }$U\in \Gamma (TN)$\textsl{\ and }$\xi \in \Gamma
(Rad(TN)).$\textsl{\ } \newline
\medskip
\noindent \textbf{Proof.~} Assume that $\nabla $ is a metric connection.
Then, $\nabla _{U}\xi \in \Gamma (Rad(TN)),$ for $U\in \Gamma (TN)$ and $\xi
\in \Gamma (Rad(TN)).$ Thus, using (\ref{eq:2.6}) and (\ref{eq:2.1c}), we
derive 
\begin{equation*}
\frac{1}{q}\breve{g}(\bar{\nabla}_{U}{\breve{J}\xi },\breve{J}U)-\frac{p}{q}%
\breve{g}(\bar{\nabla}_{U}\xi ,\breve{J}Z)=0,
\end{equation*}%
for $Z\in \Gamma (S(TN)).$

Then if we use (\ref{eq:2.6}) again and (\ref{eq:2.12}) in the equation
above, we get 
\begin{equation*}
\frac{1}{q}\breve{g}(A_{\breve{J}\xi }^{\ast }U,\breve{J}Z)-\frac{p}{q}%
\breve{g}(A_{\xi }^{\ast }U,\breve{J}Z)=0,
\end{equation*}%
which completes the proof. The converse of the assertion is obvious. \newline
\medskip
\noindent \textbf{Theorem~3.4.~}\textsl{Let }$N$\textsl{\ be an invariant
lightlike submanifold of a metallic semi-Riemannian manifold }$\breve{N}.$%
\textsl{\ Then, the radical distribution }$Rad(TN)$\textsl{\ defines a
totally geodesic foliation on }$N$\textsl{\ if and only if }%
\begin{equation*}
h^{l}(\xi ,\breve{J}U)=ph^{l}(\xi ,U),
\end{equation*}%
\textsl{for all }$U\in \Gamma (S(TN))$\textsl{\ and }$\xi \in \Gamma
(Rad(TN)).$

\medskip

\noindent \textbf{Proof.~} We assume that $Rad(TN)$ defines a totally
geodesic foliation on $N.$ That is, for $\xi ,\xi _{1}\in \Gamma (Rad(TN)),$ 
$\nabla _{\xi }{\xi _{1}}\in \Gamma (Rad(TN)).$ Since $\bar{\nabla}$ is a
metric connection, one can easily see that 
\begin{equation*}
g(\nabla _{\xi }{\xi _{1}},U)=\breve{g}(\bar{\nabla}_{\xi }{\xi _{1}},U)=%
\breve{g}(\xi _{1},\nabla _{\xi }U)=0.
\end{equation*}%
Using (\ref{eq:2.1c}), (\ref{eq:2.6}) and (\ref{eq:2.11}), we have 
\begin{equation*}
\breve{g}(\breve{J}\xi ,h^{l}(\xi _{1},\breve{J}U))-p\breve{g}(\breve{J}\xi
,h^{l}(\xi _{1},U))=0,
\end{equation*}%
and assertion is proved. The proof of the converse part can be made similarly. \newline
\medskip
\noindent \textbf{Theorem~3.5.~}\textsl{Let }$N$\textsl{\ be an invariant
lightlike submanifold of a metallic semi-Riemannian manifold }$\breve{N}.$%
\textsl{\ Then, the screen distribution }$S(TN)$\textsl{\ defines a totally
geodesic foliation on }$N$\textsl{\ if and only if }%
\begin{equation*}
{h}^{\ast }{(U,\breve{J}V)=ph}^{\ast }{(U,V)}
\end{equation*}%
\textsl{for }$U,V\in \Gamma (S(TN)).$
\medskip
\noindent \textbf{Proof.~} $S(TN)$ defines a totally geodesic foliation on $%
N $ if and only if $\nabla _{U}V\in \Gamma (S(TN))$, for $U,V\in \Gamma
(S(TN)).$ If we consider $\bar{\nabla}$ is a metric connection, we get 
\begin{equation*}
\breve{g}(\nabla _{U}V,N)=\breve{g}(\bar{\nabla}_{U}V,N)=0,
\end{equation*}%
for $N\in \Gamma (ltr(TN)).$

Using (\ref{eq:2.1c}), (\ref{eq:2.6}) and (\ref{eq:2.11}), we obtain 
\begin{eqnarray*}
&&\breve{g}(\nabla _{U}^{\ast }{\breve{J}V}+h^{\ast }(U,\breve{J}V),\breve{J}%
N) \\
&&-p.\breve{g}(\nabla _{U}^{\ast }V+h^{\ast }(U,V),\breve{J}N) \\
&=&0,
\end{eqnarray*}%
which completes the proof. The converse proof is obvious. \newline

\section*{\textbf{4.\, Screen Semi-Invariant Lightlike Submanifolds of
Metallic Semi-Riemannian Manifolds}}

\setcounter{equation}{0} \renewcommand{\theequation}{4.\arabic{equation}}

\noindent \textbf{Definition~2.~}\textsl{Let $(\breve{N},\breve{J},\breve{g}%
) $ be a metallic semi-Riemannian manifold and $(N,g)$ be a lightlike
submanifold of $\breve{N}.$ Then, we say that $N$ is a screen semi-invariant
lightlike submanifold of $\breve{N},$ if the following conditions are
satisfied:} 
\begin{equation}
\breve{J}(Rad(TN))\subseteq S(TN)  \label{eq:4.1}
\end{equation}%
and 
\begin{equation}
\breve{J}(ltr(TN))\subseteq S(TN).  \label{eq:4.2}
\end{equation}

From definition above for a screen semi-invariant lightlike submanifold of a
metallic semi-Riemannian manifold, we can define a non-degenerate
distribution $L_{0}$ such that $S(TN)$ is decomposed as: 
\begin{equation}
S(TN)=L_{0}\bot {L_{1}\oplus L_{2}},  \label{eq:4.3}
\end{equation}%
where $L_{1}=\breve{J}(Rad(TN))$ and $L_{2}=\breve{J}(ltr(TN)).$ \newline

\noindent \textbf{Proposition~4.1~}\textsl{\ Let }$(N,g)$\textsl{\ be a
screen semi-invariant lightlike submanifold of a metallic semi-Riemannian
manifold }$(\breve{N},\breve{J},\breve{g}).$\textsl{\ Then, }$\breve{J}L_{0}$%
\textsl{\ is invariant with respect to }$\breve{J}.$ \newline
\noindent \textbf{Proof.~} Let consider any vector field $U$ of $\Gamma
(L_{0}).$ From (\ref{eq:4.3}) and (\ref{eq:2.1b}), we derive 
\begin{equation}
\breve{g}(\breve{J}U,\xi )=0\,\,\ and\,\,\ \breve{g}(\breve{J}U,N)=0,
\label{eq:4.4}
\end{equation}%
for $\xi \in \Gamma (Rad(TN))$ and $N\in \Gamma (ltr(TN)).$ That is, $\breve{%
J}U\notin \Gamma (ltr(TN)\cup Rad(TN)).$

Similarly, using (\ref{eq:4.3}) and (\ref{eq:2.1c}), we have 
\begin{equation}
\breve{g}(\breve{J}U,\breve{J}\xi )=0\,\,\ and\,\,\ \breve{g}(\breve{J}U,%
\breve{J}N)=0.  \label{eq:4.5}
\end{equation}%
That is, $\breve{J}U\notin \Gamma (\breve{J}(ltr(TN))\cup \breve{J}%
(Rad(TN))) $ and proof is completed. \newline

\bigskip

Thus, $TN$ can be written as: 
\begin{equation}
TN={L_{1}\oplus L_{2}}\bot L_{0}\bot Rad(TN).  \label{eq:4.6}
\end{equation}%
If we denote the invariant distribution of $TN$ by $L$ such as 
\begin{equation}
L=L_{0}\bot Rad(TN)\bot \breve{J}(Rad(TN)),  \label{eq:4.7}
\end{equation}%
then, (\ref{eq:4.6}) is reduced to 
\begin{equation}
TN=L\oplus L_{2}.  \label{eq:4.8}
\end{equation}

\noindent \textbf{Proposition~4.2~} \textsl{Let }$(N,g)$\textsl{\ be a
screen semi-invariant lightlike submanifold of a metallic semi-Riemannian
manifold }$(\breve{N},\breve{J},\breve{g}).$\textsl{\ Then, }$S(TN^{\bot })$%
\textsl{\ is invariant with respect to }$\breve{J}.$

Thus, from (\ref{eq:2.d}) and (\ref{eq:4.8}), we have the following
decomposition: 
\begin{equation}
T\breve{N}={L_{1}\oplus L_{2}}\bot L_{0}\bot \{Rad(TN)\oplus ltr(TN)\}\bot
S(TN^{\bot }).  \label{eq:4.9}
\end{equation}

\noindent \textbf{Example~2.~} \textsl{Let $\breve{N}=R_{2}^{5}$ be a
metallic semi-Riemannian manifold of signature $(-,+,-,+,+)$ and metallic
structure $\breve{J}$ is defined as 
\begin{equation*}
\breve{J}(x_{1},x_{2},x_{3},x_{4},x_{5})=((p-\sigma )x_{1},\sigma
x_{2},\sigma x_{3},\sigma x_{4},(p-\sigma )x_{5}).
\end{equation*}%
Let $N$ be a submanifold of $(R_{2}^{5},\breve{J},\breve{g})$ given by }%
\begin{equation*}
{x}_{5}{=x}_{1}{+\sigma x}_{2}{+\sigma x}%
_{3.}
\end{equation*}%
\textsl{\ Then we get 
\begin{equation*}
U_{1}={\partial \,x_{1}}+{\partial \,x_{5}}\,,\,\,\,\ U_{2}={\partial \,x_{2}%
}+\sigma {\partial \,x_{5}}\,,\,\,\,\ U_{3}={\partial \,x_{3}}+\sigma {%
\partial \,x_{5}}\,,\,\,\,\ U_{4}={\partial \,x_{4}}.
\end{equation*}%
If we choose 
\begin{equation*}
\xi =\sigma U_{1}+U_{2}-U_{3}=(\sigma {\partial \,x_{1}}+{\partial \,x_{2}}-{%
\partial \,x_{3}}-\sigma {\partial \,x_{5}})
\end{equation*}%
and 
\begin{equation*}
W_{1}=U_{4}\,,\,\,\,\ W_{2}=U_{1}+U_{2}+U_{3}
\end{equation*}%
then $S(TN)=Sp\{W_{1},W_{2}\}$ and $Rad(TN)=Sp\{\xi \}.$ By direct
calculations we get }%
\begin{equation*}
{N=\frac{1}{2}\{(-\sigma {\partial \,x_{1}}+{\partial \,x_{2}}+{%
\partial \,x_{3}}+\sigma {\partial \,x_{5}}).}
\end{equation*}%
\textsl{\ Thus, $N$ is a screen semi-invariant lightlike hypersurface of $%
\breve{N}.$ }

Let us denote the projection morphisms on $L$ and $L_{2}$ by $B$ and $R,$
respectively. Then $U\in \Gamma (TN)$ can be written as: 
\begin{equation}
U=BU+RU,  \label{eq:4.10}
\end{equation}%
where $BU\in \Gamma (L)$ and $RU\in \Gamma (L_{2}).$

Appliying $\breve{J}$ to (\ref{eq:4.10}) we derive 
\begin{equation}
\breve{J}U=\breve{J}BU+\breve{J}RU.  \label{eq:4.11}
\end{equation}%
If we denote $\breve{J}BU$ and $\breve{J}RU$ by $S_{1}U$ and $R_{1}U,$
respectively, then we can rewrite (\ref{eq:4.11}) as 
\begin{equation}
\breve{J}U=S_{1}U+R_{1}U,  \label{eq:4.12}
\end{equation}%
where $S_{1}U\in \Gamma (L)$ and $R_{1}U\in \Gamma (ltr(TN)).$ \newline

Let $N$ be a screen semi-invariant lightlike submanifold of a metallic
semi-Riemannian manifold $\breve{N}.$ Using (\ref{eq:3.5}), (\ref{eq:2.6})
and (\ref{eq:2.11}), $\forall U,V\in \Gamma (TN),$ we obtain 
\begin{eqnarray}
\breve{J}\nabla _{U}V+\breve{J}h^{l}(U,V)+\breve{J}h^{s}(U,V) &=&\nabla
_{U}^{\ast }\breve{J}V+h^{\ast }(U,\breve{J}V)  \notag \\
&+&h^{l}(U,\breve{J}V)+h^{s}(U,\breve{J}V).  \label{eq:4.13}
\end{eqnarray}
Considering the tangential, lightlike transversal and screen transversal
parts of (\ref{eq:4.13}), we give following.

\medskip

\noindent\textbf{Proposition~4.3}~\textsl{Let }$N$\textsl{\ be a screen
semi-invariant lightlike submanifold of a metallic semi-Riemannian manifold }%
$\breve{N}.$\textsl{\ Then we have }%
\begin{equation}
\breve{J}\nabla _{U}V=\nabla _{U}^{\ast }\breve{J}V+h^{\ast }(U,\breve{J}V)-%
\breve{J}h^{l}(U,V),  \label{eq:4.14}
\end{equation}%
\begin{equation}
h^{l}(U,\breve{J}V)=0,  \label{eq:4.15}
\end{equation}%
\begin{equation}
h^{s}(U,\breve{J}V)=\breve{J}h^{s}(U,V),  \label{eq:4.16}
\end{equation}%
\textsl{for all }$U,V\in \Gamma (TN).$

\medskip

\noindent \textbf{Theorem~4.1.}~\textsl{Let }$N$\textsl{\ be a screen
semi-invariant lightlike submanifold of a metallic semi-Riemannian manifold }%
$\breve{N}.$\textsl{\ Then, the invariant distribution }$L$\textsl{\ is
integrable if and only if }%
\begin{equation*}
h^{l}(\breve{J}W,\breve{J}U)=ph^{l}(U,\breve{J}W)+qh^{l}(U,W)
\end{equation*}%
\textsl{for all }$U,W\in \Gamma (L).$

\medskip

\noindent \textbf{Proof.}~ $L$ is integrable iff $[U,W]\in \Gamma (L),$ for
all $U,W\in \Gamma (L).$ That is, $\breve{g}([\breve{J}W,U],\breve{J}\xi
)=0, $ for $\xi \in \Gamma (Rad(TN)).$ Then, using (\ref{eq:2.1b}), (\ref%
{eq:3.5}) and (\ref{eq:2.1c}), we get 
\begin{equation*}
\breve{g}(\bar{\nabla}_{\breve{J}W}{\breve{J}U},\xi )-p\breve{g}(\bar{\nabla}%
_{U}W,\breve{J}\xi )-q\breve{g}(\bar{\nabla}_{U}W,\xi )=0.
\end{equation*}%
Finally, if we use (\ref{eq:2.6}) in the last equation, the proof is
completed. \newline

\medskip

\noindent \textbf{Theorem~4.2.}~\textsl{Let }$N$\textsl{\ be a screen
semi-invariant lightlike submanifold of a metallic semi-Riemannian manifold }%
$\breve{N}.$\textsl{\ Then, the radical distribution }$Rad(TN)$\textsl{\ is
integrable if and only if }%
\begin{equation}
\nabla _{U}^{\ast }\breve{J}W-\nabla _{W}^{\ast }\breve{J}U=p(A_{U}^{\ast
}W-A_{W}^{\ast }U)  \label{eq:4.17}
\end{equation}%
\textsl{or }%
\begin{equation}
\breve{J}\nabla _{U}^{\ast }\breve{J}W-\breve{J}\nabla _{W}^{\ast }\breve{J}%
U=p(\nabla _{U}^{\ast }\breve{J}W-\nabla _{W}^{\ast }\breve{J}U)
\label{eq:4.18}
\end{equation}%
\textsl{for any }$U,W\in \Gamma (Rad(TN)).$

\medskip

\noindent \textbf{Proof.}~ We assume that $Rad(TN)$ is integrable. Then, $%
g([U,W],Z)=0$ for for all $U,W\in \Gamma (Rad(TN))\,,\,\ Z\in \Gamma
(S(TN)). $ If we use (\ref{eq:2.1c}), we have 
\begin{eqnarray}
0 &=&\breve{g}(\bar{\nabla}_{U}{\breve{J}W},\breve{J}Z)-p\breve{g}(\bar{%
\nabla}_{U}W,\breve{J}Z)  \notag \\
&-&\breve{g}(\bar{\nabla}_{W}{\breve{J}U},\breve{J}Z)+p\breve{g}(\bar{\nabla}%
_{W}U,\breve{J}Z)  \label{eq:4.19}
\end{eqnarray}%
Then, using (\ref{eq:2.6}) and (\ref{eq:2.11}), we obtain 
\begin{equation*}
\breve{g}(\nabla _{U}^{\ast }{\breve{J}W}+pA_{W}^{\ast }U-\nabla _{W}^{\ast }%
{\breve{J}U}-pA_{U}^{\ast }W,\breve{J}Z)=0,
\end{equation*}%
which satisfies (\ref{eq:4.17}).

On the other hand, if we use (\ref{eq:2.1b}), (\ref{eq:2.6}) and (\ref%
{eq:2.11}) in (\ref{eq:4.19}), we get 
\begin{eqnarray*}
0 &=&\breve{g}(\nabla _{U}^{\ast }{\breve{J}W}-\nabla _{W}^{\ast }{\breve{J}U%
},\breve{J}Z) \\
&+&p\breve{g}(\nabla _{W}^{\ast }{\breve{J}U}+h^{\ast }(W,\breve{J}U)-\breve{%
J}h^{l}(W,U)-\nabla _{U}^{\ast }{\breve{J}W}+h^{\ast }(U,\breve{J}W)-\breve{J%
}h^{l}(U,W),Z).
\end{eqnarray*}%
Since $h^{l}$ is symmetric and from (\ref{eq:4.14}), we have 
\begin{equation*}
\breve{g}(\breve{J}\nabla _{U}^{\ast }{\breve{J}W}-\nabla _{W}^{\ast }{%
\breve{J}U}+p\nabla _{W}^{\ast }{\breve{J}U}-p\nabla _{U}^{\ast }{\breve{J}W}%
,Z)=0
\end{equation*}%
which satisfies (\ref{eq:4.18}) and the proof is completed. \newline
\medskip
\noindent \textbf{Theorem~4.3.}~\textsl{Let }$N$\textsl{\ be a screen
semi-invariant lightlike submanifold of a metallic semi-Riemannian manifold }%
$\breve{N}.$\textsl{\ Then, the screen distribution }$S(TN)$\textsl{\ is
integrable if and only if }%
\begin{equation}
\nabla _{U}^{\ast }\breve{J}W-\nabla _{W}^{\ast }\breve{J}U=p(\nabla
_{U}^{\ast }W-\nabla _{W}^{\ast }U)  \label{eq:4.20}
\end{equation}%
\textsl{or }%
\begin{equation}
\nabla _{U}^{\ast }\breve{J}W=\nabla _{W}^{\ast }\breve{J}U)  \label{eq:4.21}
\end{equation}%
\textsl{for any }$U,W\in \Gamma (S(TN)).$
\medskip
\noindent \textbf{Proof.}~ We know that $Rad(TN)$ is integrable iff $\breve{%
g}([U,W],N)=0,$ for for all $U,W\in \Gamma (S(TN))\,,\,\ N\in \Gamma
(ltr(TN)).$ Using (\ref{eq:2.1c}) and (\ref{eq:3.5}), we get 
\begin{equation}
\breve{g}(\bar{\nabla}_{U}{\breve{J}W}-\bar{\nabla}_{W}{\breve{J}U},\breve{J}%
N)-p\breve{g}(\bar{\nabla}_{U}W-\bar{\nabla}_{W}U,\breve{J}N)=0.
\label{eq:4.22}
\end{equation}%
If we use (\ref{eq:2.6}) and (\ref{eq:2.11}) in (\ref{eq:4.22}), 
\begin{equation*}
\breve{g}(\nabla _{U}^{\ast }{\breve{J}W}-\nabla _{W}^{\ast }{\breve{J}U}%
-p\nabla _{U}^{\ast }W+p\nabla _{W}^{\ast }U,\breve{J}N)=0
\end{equation*}%
is obtained and (\ref{eq:4.20}) is satisfies.

On the other hand, using (\ref{eq:2.1b}), (\ref{eq:2.6}) and (\ref{eq:2.11})
in (\ref{eq:4.22}), we get 
\begin{equation*}
\breve{g}(\nabla _{U}^{\ast }{\breve{J}W}-\nabla _{W}^{\ast }{\breve{J}U},%
\breve{J}N)=0.
\end{equation*}%
Thus, (\ref{eq:4.21}) is satisfied and the proof is completed. \newline
\medskip
\noindent \textbf{Theorem~4.4.}~\textsl{Let }$N$\textsl{\ be a screen
semi-invariant lightlike submanifold of a metallic semi-Riemannian manifold }%
$\breve{N}.$\textsl{\ Then, the induced connection }$\nabla $\textsl{\ on }$%
N $\textsl{\ is a metric connection if and only if one of the followings is
satisfied: }%
\begin{equation}
\nabla _{U}^{\ast }{\breve{J}\xi }=-pA_{\xi }^{\ast }U  \label{eq:4.23}
\end{equation}%
\textsl{or }%
\begin{equation}
qA_{\xi }^{\ast }U=0  \label{eq:4.24}
\end{equation}%
\textsl{for all }$U\in \Gamma (TN)$\textsl{\ and }$\xi \in \Gamma (Rad(TN)).$%

\medskip
\noindent \textbf{Proof.}~ $\nabla $ is a metric connection iff $\nabla
_{U}\xi \in \Gamma (Rad(TN))$ that is, $g(\nabla _{U}\xi ,Z)=0,$ for $U\in
\Gamma (TN),$ $\xi \in \Gamma (Rad(TN))$ and $Z\in \Gamma (S(TN)).$ If we
consider (\ref{eq:2.6}), we have 
\begin{equation}
\breve{g}(\bar{\nabla}_{U}\xi ,Z)=0.  \label{eq:4.25}
\end{equation}%
Then, using (\ref{eq:2.1c}), (\ref{eq:3.5}), (\ref{eq:2.6}) and (\ref%
{eq:2.11}), we obtain 
\begin{equation*}
\breve{g}(\nabla _{U}^{\ast }{\breve{J}\xi }+pA_{\xi }^{\ast }U,\breve{J}Z)=0
\end{equation*}%
which satisfies (\ref{eq:4.23}).

Now, using (\ref{eq:2.1c}) in (\ref{eq:4.25}), we get 
\begin{equation*}
\breve{g}(\breve{J}\bar{\nabla}_{U}{\breve{J}\xi }-p\breve{J}\bar{\nabla}%
_{U}\xi ,Z)=0.
\end{equation*}%
Using (\ref{eq:2.6}), (\ref{eq:2.11}) and (\ref{eq:2.12}), we derive 
\begin{eqnarray*}
0 &=&\breve{g}(\breve{J}\nabla _{U}^{\ast }{\breve{J}\xi }+\breve{J}h^{\ast
}(U,\breve{J}\xi )+\breve{J}h^{l}(U,\breve{J}\xi )+\breve{J}h^{s}(U,\breve{J}%
\xi ) \\
&+&p\breve{J}A_{\xi }^{\ast }U-p\breve{J}\nabla _{U}^{\ast t}\xi -p\breve{J}%
h^{l}(U,\xi )-p\breve{J}h^{s}(U,\xi ),Z).
\end{eqnarray*}%
Finally, if we use (\ref{eq:2.1c}), (\ref{eq:4.15}) and (\ref{eq:4.16}) in
the last equation, 
\begin{equation*}
\breve{g}(-qA_{\xi }^{\ast }U+q\nabla _{U}^{\ast t}\xi ,Z)=0
\end{equation*}%
is obtained and (\ref{eq:4.24}) is satisfied. \newline
\medskip
\noindent \textbf{Theorem~4.5.}~\textsl{Let }$N$\textsl{\ be a screen
semi-invariant lightlike submanifold of a metallic semi-Riemannian manifold }%
$\breve{N}.$\textsl{\ Then, the radical distribution }$Rad(TN)$\textsl{\
defines a totally geodesic foliation on }$N$\textsl{\ if and only if }%
\begin{equation*}
\nabla _{\xi }^{\ast }{\breve{J}U}=p\nabla _{\xi }^{\ast }U,
\end{equation*}%
\textsl{for all }$U\in \Gamma (S(TN))$\textsl{\ and }$\xi \in \Gamma
(Rad(TN)).$
\medskip
\noindent \textbf{Proof.}~ We assume that $Rad(TN)$ defines a totally
geodesic foliation on $N.$ That is, for $\xi ,\xi _{1}\in \Gamma (Rad(TN)),$ 
$\nabla _{\xi }\xi _{1}\in \Gamma (Rad(TN)).$ Since $\bar{\nabla}$ is a
metric connection, we can write 
\begin{equation*}
g(\nabla _{\xi }{\xi _{1}},U)=\breve{g}(\bar{\nabla}_{\xi }{\xi _{1}},U)=%
\breve{g}(\xi _{1},\bar{\nabla}_{\xi }U)=0.
\end{equation*}%
Using (\ref{eq:2.1c}), (\ref{eq:2.6}) and (\ref{eq:2.11}), we have 
\begin{equation*}
\breve{g}(\breve{J}\xi ,\nabla _{\xi }^{\ast }{\breve{J}U}-p\nabla _{\xi
}^{\ast }U)=0
\end{equation*}%
and assertion is proved.

\medskip

\noindent \textbf{Theorem~4.6.}~\textsl{Let }$N$\textsl{\ be a screen
semi-invariant lightlike submanifold of a metallic semi-Riemannian manifold }%
$\breve{N}.$\textsl{\ Then, the screen distribution }$S(TN)$\textsl{\
defines a totally geodesic foliation on }$N$\textsl{\ if and only if }%
\begin{equation*}
\nabla _{U}^{\ast }{\breve{J}V}=p.\nabla _{U}^{\ast }V
\end{equation*}%
\textsl{for all }$U,V\in \Gamma (S(TN)).$
\medskip
\noindent \textbf{Proof.}~ $S(TN)$ defines a totally geodesic foliation on $%
N $ iff for $U,V\in \Gamma (S(TN)),$ $\nabla _{U}V\in \Gamma (S(TN)).$ If we
consider $\bar{\nabla}$ is a metric connection, we get 
\begin{equation*}
\breve{g}(\nabla _{U}V,N)=\breve{g}(\bar{\nabla}_{U}V,N)=0,
\end{equation*}%
for $N\in \Gamma (ltr(TN)).$ Using (\ref{eq:2.1c}), (\ref{eq:2.6}) and (\ref%
{eq:2.11}), 
\begin{equation*}
\breve{g}(\nabla _{U}^{\ast }{\breve{J}V},\breve{J}N)-p.\breve{g}(\nabla
_{U}^{\ast }V,\breve{J}N)=0,
\end{equation*}%
which completes the proof.

\bigskip

{\small \textbf{Bilal Eftal Acet} \newline
Department of Mathematics \newline
Ad\i yaman University, 02040, Ad\i yaman, Turkey \newline
e-mail: eacet@adiyaman.edu.tr }

{\small {\ \textbf{Feyza Esra Erdo\u{g}an} \newline
Faculty Science, Department of Mathematics \newline
Ege University, \.{I}zmir, Turkey \newline
e-mail: feyza.esra.erdogan@ege.edu.tr }}

{\small \smallskip }

{\small \textbf{Selcen Y\"{u}ksel Perkta\c{s}} \newline
Department of Mathematics, Faculty of Arts and Sciences \newline
Ad\i yaman University, 02040, Ad\i yaman, Turkey \newline
e-mail: sperktas@adiyaman.edu.tr }

{\small \smallskip }

\end{document}